
\documentstyle{amsppt}
\magnification=1200
\def\ssbull{\raise.3ex\hbox{${\scriptscriptstyle\bullet}$}}
\def\scirc{{\scriptstyle{\circ}}}
\def\Sym{\hbox{{\rm Sym}}}
\def\cHom{\hbox{${\Cal H}om $}}

\def\Gr{\text{{\rm Gr}}}
\def\Spec{\hbox{{\rm Spec}}}
\def\Im{\hbox{{\rm Im}}}

\def\Ker{\hbox{{\rm Ker}}}
\def\Hom{\hbox{{\rm Hom}}}

\def\DR{\hbox{{\rm DR}}}
\def\Diff{\text{{\rm Diff}}}
\def\supp{\hbox{{\rm supp}}\,}
\def\Supp{\hbox{{\rm Supp}}\,}
\def\codim{\text{\rm codim}}
\def\id{\text{\rm id}}
\def\qcoh{\text{\rm qcoh}}
\def\red{\text{\rm red}}
\def\an{\text{\rm an}}

\def\prim{\text{\rm prim}}
\def\simto{\buildrel\sim\over\to}
\def\SameAuthor{\vrule height3pt depth-2.5pt width1cm}

\topmatter
\title Dwork Cohomology and Algebraic $ {\Cal D} $-Modules
\endtitle
\author A.~Dimca,  F.~Maaref, C.~Sabbah and M.~Saito\endauthor
\abstract Using local cohomology and algebraic 
$ {\Cal D} $-Modules, we generalize a comparison theorem 
between relative de Rham cohomology and Dwork cohomology due to 
N.~Katz, A.~Adolphson and S.~Sperber.
\endabstract
\endtopmatter
\tolerance=1000
\baselineskip=12pt
\vsize=23truecm

\document

\centerline{\bf Introduction}

\bigskip\noindent
Let
$ p : X \rightarrow  S $ be a smooth affine morphism of smooth complex 
algebraic varieties, and
$ f = (f_{1}, \dots , f_{r}) : X \rightarrow {\Bbb A}^{r} $ a morphism to 
the affine space
$ {\Bbb A}^{r} $ with coordinates
$ (y_{1}, \dots , y_{r}) $.
We define
$ V = X \times  {\Bbb A}^{r}, F = \sum _{j} f_{j} (x) y_{j} $ and
$ Y = f^{-1}(0)_{\red} $ so that we have a diagram
$$
\CD
Y @>i>> X   @<\pi <<  V  @>F>>  {\Bbb A}^{1}
\\
@.  @VVpV 
\\
 @.  S
\endCD
$$
where
$ i $ and
$ \pi  $ denote natural morphisms.
(More precisely,
$ f $ is a morphism to the dual space
$ ({\Bbb A}^{r})^{\vee } $ of
$ {\Bbb A}^{r} $ which is identified with
$ {\Bbb A}^{r} $,
and
$ F $ is the pull-back of the canonical bilinear form on
$ ({\Bbb A}^{r})^{\vee } \times {\Bbb A}^{r} $ by
$ f \times \id $.)

Let
$ {\Cal E} $ be a locally free
$ {\Cal O}_{X} $-Module of finite rank with an integrable connection
$ \nabla  $.
We first assume that
$ Y $ is smooth over
$ S $,
and let
$ d =  \codim _{X}Y $.
Then the restriction
$ i^{*}{\Cal E} $ has the induced connection which is also denoted by
$ \nabla  $.
We define the twisted connection
$ \nabla _{F} $ on the pull-back
$ \pi ^{*}{\Cal E} $ by
$ \nabla _{F}e = \nabla e + dF\otimes e $.
Since
$ p\scirc i $ and
$ p\scirc \pi  $ are smooth, we have the relative de Rham 
cohomologies
$ R^{j}(p\scirc i)_{*}\DR_{Y/S}(i^{*}{\Cal E},\nabla ) $,
$ R^{k}(p\scirc \pi )_{*}\DR_{V/S}(\pi ^{*}{\Cal E},\nabla _{F}) $ which 
are 
$ {\Cal O}_{S} $-Modules with the Gauss-Manin connection as defined in 
[8].

Assume further that
$ f^{-1}(0) $ is reduced and is a complete intersection (i.e.,
$ d = r) $.
Then A.~Adolphson and S.~Sperber [1] proved the following canonical 
isomorphism of
$ {\Cal O}_{S} $-Modules with an integrable connection :
$$
R^{j}(p\scirc i)_{*}\DR_{Y/S}(i^{*}{\Cal E},\nabla ) = 
R^{j+2d}(p\scirc \pi )_{*}\DR_{V/S}(\pi ^{*}{\Cal E},\nabla _{F}).
\tag 0.1
$$
This is a generalization of [7] in the case
$ r = 1 $.
The isomorphism means that the relative de Rham cohomology with the 
Gauss-Manin connection on the left-hand side can be calculated by the 
right-hand side which is called the {\it Dwork cohomology} associated 
with
$ F $.
For example, if
$ X = {\Bbb A}^{n} \times  S $ with
$ S =  \Spec \, A $,
$ p : X \rightarrow  S $ is the natural projection, and
$ ({\Cal E},\nabla ) = ({\Cal O}_{X},d) $,
then the Dwork cohomology is identified with the cohomology of the 
Koszul complex
$$
K(A[x,y]; \partial /\partial x_{j} + \sum _{i} y_{i}\partial f_{i}/
\partial x_{j} \,(1\le  j \le  n), 
\partial /\partial y_{i} + f_{i} \,(1\le  i \le  r)),
$$
where
$ A[x,y] = A[x_{1}, \dots , x_{n}, y_{1}, \dots , y_{r}] $ is the affine 
coordinate ring of
$ S \times  {\Bbb A}^{n} \times  {\Bbb A}^{r} $.

In this paper, we give a generalization of (0.1) which holds without any 
assumption on
$ p : X \rightarrow  S $ and
$ Y $ (but
$ X $ and
$ S $ remain smooth).
The reader would, however, notice soon that (0.1) is useless if
$ p $ is proper, because the restrictions of
$ f_{i} $ to the fibers of
$ p $ are constant in this case.
So
$ V = X \times  {\Bbb A}^{r} $ should be replaced by an algebraic vector 
bundle
$ \pi  : V \rightarrow  X $ with rank
$ r $,
and
$ f = (f_{1}, \dots , f_{r}) $ by a section
$ s : X \rightarrow  V^{\vee } $ of the dual vector bundle
$ \pi ^{\vee } : V^{\vee } \rightarrow  X $ of
$ V $.
Then
$ s $ is naturally identified with a function on
$ V $ such that the restriction to each fiber of
$ V $ is a linear function, and this function plays the role of
$ F $.
See (2.1).
The subspace
$ Y $ is now defined as
$ s^{-1}(0)_{\red} $,
where
$ 0 $ denotes the zero section.
If
$ V $ is trivial (by shrinking
$ X $ if necessary), we recover the previous situation.

The main idea is to use the formalism of algebraic
$ {\Cal D} $-Modules rather than that of relative de Rham cohomology 
and Gauss-Manin connection.
The reason is that if
$ Y $ is singular, the restriction
$ i^{*}({\Cal E},\nabla ) $ should be replaced with 
$ \bold{R}\Gamma _{Y}({\Cal E},\nabla )[d] $,
but this is useful only if it is defined in the derived category of algebraic
$ {\Cal D}_{X} $-Modules.
So it is natural to replace also
$ ({\Cal E},\nabla ) $ with a bounded complex of quasi-coherent left
$ {\Cal D}_{X} $-Modules
$ M^{\ssbull } $.
Then the usual pull-back
$ \pi ^{*}M^{\ssbull } $ has naturally a structure of (complex of) 
quasi-coherent left
$ {\Cal D}_{V} $-Modules.
We will denote by
$ (\pi ^{*}M^{\ssbull })_{F} $ the complex of
$ {\Cal D}_{V} $-Modules with twisted action of
$ {\Cal D}_{V} $ as above.

Another merit of using
$ {\Cal D} $-Modules is that we can reduce the assertion to the case
$ p  $ is the identity due to the compatibility of the direct image functor 
with the composition of morphisms.

\proclaim{0.2.~Theorem}
We have a canonical isomorphism in the derived category of left
$ {\Cal D}_{X} $-Modules
$$
\bold{R}\Gamma _{Y}M^{\ssbull }[r] = \pi _{+}((\pi ^{*}M^{\ssbull })_{F}).
$$
\endproclaim

Here
$ \pi _{+} $ denotes the direct image of algebraic left
$ {\Cal D} $-Modules (see e.g.~[2]).
If
$ M^{\ssbull } = ({\Cal E},\nabla ) $ and
$ Y $ is smooth with codimension
$ d $,
then
$ \bold{R}\Gamma _{Y}M^{\ssbull }[d] = i_{+}(i^{*}M^{\ssbull }) $.

In order to compare Theorem (0.2) with the result of A.~Adolphson and 
S.~Sperber (0.1), we apply the direct image functor
$ p_{+} $ to both sides of (0.2), and use the isomorphism
$ (p\scirc \pi )_{+} = p_{+}\pi _{+} $.
However, we need to know that the cohomology of the direct image of a
$ {\Cal D} $-Module under a smooth morphism is naturally isomorphic to 
the relative de Rham complex with the Gauss-Manin connection as 
defined by Katz-Oda [8].
We give in (1.4) a proof of this assertion, which is apparently widely 
accepted, but for which no complete proof seems to exist in the 
literature.
From (0.2) and (1.4) we can deduce the following (see Remark (2.5)(i)) :

\proclaim{0.3.~Corollary}
The isomorphism (0.1) remains true for any smooth morphism
$ p $ and any section
$ s $ of an algebraic vector bundle
$ \pi  : V \rightarrow  X $ such that
$ Y \,(= s^{-1}(0)_{\red}) $ is smooth over
$ S $.
\endproclaim

Note that
$ s^{-1}(0) $ may be nonreduced, and the codimension of
$ Y $ may be different from the rank of
$ V $ (in particular,
$ Y $ is not necessarily a complete intersection).
It does not seem easy to prove the assertion under such an assumption 
without using the theory of
$ {\Cal D} $-Modules.

We give two proofs of (0.2).
The first proof (2.2) reduces the assertion to a lemma on Fourier 
transformation (Lemma (2.3)).
In the case where
$ V = X \times  {\Bbb A}^{r}, f = (f_{1}, \dots , f_{r}) : X \rightarrow  
({\Bbb A}^{r})^{\vee } $ is an isomorphism, and
$ M^{\ssbull } $ is a quasi-coherent
$ {\Cal D} $-Module
$ M $,
it is well-known to specialists that the Fourier transform
$ \widehat{M} $ of
$ M $ is isomorphic to
$ \pi '_{+}(\pi ^{*}M)_{F} $,
where
$ \pi ' : X \times  {\Bbb A}^{r} \rightarrow  {\Bbb A}^{r} $ is the second 
projection.
So Theorem (0.2) in this case means that the algebraic de Rham 
cohomology of
$ \widehat{M} $ is isomorphic to the cohomology of the restriction of
$ M $ to the origin of
$ {\Bbb A}^{r} $ as
$ {\Cal D} $-Module (or equivalently as
$ {\Cal O} $-Module).

The second proof (2.4) reduces the assertion to a calculation of local 
cohomology, and is used in the proof of (0.4).
(The two proofs are both contained in this paper, because it is not yet 
shown that the obtained isomorphisms coincide.)

Assume now
$ X = {\Bbb P}_{S}^{n} $ with
$ p : X \rightarrow  S $ the natural projection, and
$ Y $ is a hypersurface defined by a homogeneous polynomial
$ F $ of degree
$ m $ with coefficients in
$ \Gamma (S,{\Cal O}_{S}) $.
Then
$ F $ is identified with a global section of
$ {\Cal O}_{X}(m) \,(= {\Cal O}_{X}(V^{\vee })) $.
Let
$ (x_{0}, \dots , x_{n}) $ be the natural coordinate of
$ {\Bbb A}^{n+1} $ so that
$ \Gamma ({\Bbb A}^{n+1}, {\Cal O}) = {\Bbb C}[x] 
\,(:= {\Bbb C}[x_{0}, \dots , x_{n}]) $.
Let
$$
\Omega ^{i}({\Bbb A}_{S}^{n+1}/S)^{(m)} = \sum \Sb |I|=i \\ 
|\nu |+|I|\equiv 0 \mod m \endSb {\Cal O}_{S}x^{\nu }dx_{I},
$$
where
$ x^{\nu } = \prod _{i} {x}_{i}^{{\nu }_{i}} $ for
$ \nu  = (\nu _{0}, \dots , \nu _{n}) \in  {\Bbb Z}^{n+1}  $ and
$ dx_{I} = dx_{j_{1}}\wedge  \dots \wedge dx_{j_{i}} $ for
$ I = \{j_{1}, \dots , j_{n}\} $.
Then
$ \Omega ^{\ssbull }({\Bbb A}_{S}^{n+1}/S)^{(m)} $ is stable by
$ d + dF\wedge  $.
We call
$$
{\Cal H}^{i}(\Omega ^{\ssbull }({\Bbb A}_{S}^{n+1}/S)^{(m)}, d + dF\wedge )
$$
the primitive Dwork cohomology sheaf of
$ F $.
It has a left
$ {\Cal D}_{S} $-Module structure such that the action of a vector field
$ \xi  $ on
$ S $ is given by
$ m \mapsto  \xi m + (\xi F)m $,
where
$ \xi m $ and
$ \xi F $ are defined by using the natural action of
$ \xi  $ on
$ {\Cal O}_{S} $.

We define the primitive part of the direct image
$ {\Cal H}^{i}p_{+}\bold{R}\Gamma _{Y}{\Cal O}_{X} $ of the left
$ {\Cal D}_{X} $-Module
$ \bold{R}\Gamma _{Y}{\Cal O}_{X} $ by
$$
({\Cal H}^{i}p_{+}\bold{R}\Gamma _{Y}{\Cal O}_{X})^{\prim} = 
\Ker({\Cal H}^{i}p_{+}\bold{R}\Gamma _{Y}{\Cal O}_{X} 
\rightarrow  {\Cal H}^{i}p_{+}{\Cal O}_{X}).
$$
Note that the morphism on the right-hand side is surjective for
$ i \ne  - n $ and zero for
$ i = - n $.
See (3.1).

\proclaim{0.4.~Theorem}
We have a canonical isomorphism of left
$ {\Cal D}_{S} $-Modules
$$
({\Cal H}^{i}p_{+}\bold{R}\Gamma _{Y}{\Cal O}_{X})^{\prim} = 
{\Cal H}^{i+n}(\Omega ^{\ssbull }({\Bbb A}_{S}^{n+1}/S)^{(m)}, d + dF\wedge ).
$$
\endproclaim

So the left-hand side is calculated by a subcomplex of the Koszul 
complex for
$ \partial /\partial x_{i} + \partial F/\partial x_{i} \,(0 \le  i \le  n) $.
This theorem is closely related with results of Katz [7].
Note that Theorem (0.4) in the case
$ S = pt $ becomes
$$
\tilde{H}^{i-1}(X \backslash  Y, \bold{C}) = {H}_{Y}^{i}(X,{\Bbb C})^{\prim} 
= H^{i}(\Omega ^{\ssbull }({\Bbb A}_{S}^{n+1}/S)^{(m)}, d + dF\wedge ),
$$
and follows from [3], because
$ X \backslash  Y $ is the quotient of
$ F^{-1}(1) \subset  {\Bbb A}^{n+1} $ by the action of the geometric 
monodromy automorphism.
Theorem (0.4) is also related with Griffiths' result [5] in the case
$ Y $ is smooth over
$ S $.
See (3.4).

In \S 1, we review elementary facts from the theory of algebraic
$ {\Cal D} $-Modules [2] together with the relation with differential 
complexes, and prove the equivalence between the direct image of 
algebraic  $ {\Cal D} $-Modules and the Gauss-Manin connection of 
Katz-Oda [8] in the case of smooth morphisms.
Then we show Theorems (0.2) and (0.4) in \S 2 and \S 3.

\bigskip\noindent
{\it Acknowledgement.}
We would like to thank A. Adolphson and S. Sperber for sending us 
promptly their preprint [1].
During the preparation of this paper, the first three authors have
benefited from the financial support by INTAS program 97-1644.

\bigskip\bigskip\centerline{{\bf 1.~Review on Algebraic $ {\Cal D} $-Modules}}

\bigskip
\noindent
{\bf 1.1.}
{\it $ {\Cal D} $-Modules.} Let
$ X $ be a smooth complex algebraic variety with dimension
$ n $,
and
$ {\Cal D}_{X} $ the sheaf of linear algebraic differential operators 
with
$ F $ the filtration by the order of differential operators.
We denote by
$ M({\Cal D}_{X})^{l} $ the category of left
$ {\Cal D}_{X} $-Modules, and by
$ M_{\qcoh}({\Cal D}_{X})^{l} $ the full subcategory consisting of 
quasi-coherent left
$ {\Cal D}_{X} $-Modules.
We have
$ C^{*}({\Cal D}_{X})^{l}, K^{*}({\Cal D}_{X})^{l}, D^{*}({\Cal D}_{X})^{l} $ 
for
$ * = +, -, b $ or empty as usual [2], [12].
We denote by
$ {D}_{\qcoh}^{{}^*}({\Cal D}_{X})^{l} $ the full subcategory of
$ D^{*}({\Cal D}_{X})^{l} $ consisting of complexes with quasi-coherent 
cohomologies.
It is known that
$ {D}_{\qcoh}^{b}({\Cal D}_{X})^{l} $ is equivalent to the derived 
category consisting of bounded complexes of quasi-coherent
$ {\Cal D}_{X} $-Modules.
See [2, VI, 2.10].

We have similarly
$ M({\Cal D}_{X})^{r}, D^{*}({\Cal D}_{X})^{r}, {D}_{\qcoh}^{{}^*}
({\Cal D}_{X})^{r},  $ etc.
for right
$ {\Cal D}_{X} $-Modules.
We have the transformation between the left and right
$ {\Cal D} $-Modules by
$$
M \rightarrow  \omega _{X}\otimes _{{\Cal O}_{X}}M
\tag 1.1.1
$$
for a left
$ {\Cal D}_{X} $-Module
$ M $,
where
$ \omega _{X} $ is the dualizing sheaf (i.e.
 $ {\Omega }_{X}^{\dim X}) $.
See loc.~cit.

Note that a left
$ {\Cal D}_{X} $-Module can be viewed as an
$ {\Cal O}_{X} $-Module (not necessarily finitely generated) endowed 
with an integrable connection.
By the theory of characteristic variety, it is known that a left
$ {\Cal D}_{X} $-Module is locally free of finite type over
$ {\Cal O}_{X} $ if and only if it is coherent over
$ {\Cal O}_{X} $.

\bigskip
\noindent
{\bf 1.2.}
{\it Differential Complexes.} Let
$ L, L' $ be
$ {\Cal O}_{X} $-Modules.
Then we have right
$ {\Cal D}_{X} $-Modules
$ L\otimes _{{\Cal O}_{X}}{\Cal D}_{X}, L'\otimes _{{\Cal O}_{X}}
{\Cal D}_{X} $,
and any morphism
$$
u \in  \Hom_{{\Cal D}_{X}}(L\otimes _{{\Cal O}_{X}}{\Cal D}_{X}, 
L'\otimes _{{\Cal O}_{X}}{\Cal D}_{X}) = \Hom_{{\Cal O}_{X}}(L, 
L'\otimes _{{\Cal O}_{X}}{\Cal D}_{X})
$$
induces a  
$ {\Bbb C} $-linear morphism
$ u' : L \rightarrow  L' $ by taking the composition with the the 
projection
$ L'\otimes _{{\Cal O}_{X}}{\Cal D}_{X} \rightarrow  L' $ (defined by
$ m\otimes P \mapsto  (P1)m) $.
We see that
$ u' $ is uniquely determined by the commutative diagram
$$
\CD
L\otimes _{{\Cal O}_{X}}{\Cal D}_{X} @>u>>
L'\otimes _{{\Cal O}_{X}}{\Cal D}_{X}
\\
@VVV @VVV
\\
L  @>u' >>   L'
\endCD
\tag 1.2.1
$$
where the vertical morphisms are as above.
The obtained morphism
$$
\Hom_{{\Cal D}_{X}}(L\otimes _{{\Cal O}_{X}}{\Cal D}_{X}, 
L'\otimes _{{\Cal O}_{X}}{\Cal D}_{X}) \rightarrow  
\Hom_{{\Bbb C}}(L, L')
\tag 1.2.2
$$
is injective.
See [10, 2.2.2].
Its image is denoted by
$ \Diff _{X}(L,L') $.
It has an increasing filtration
$ F $ with
$ p $-th term defined by
$$
F_{p} \Diff _{X}(L, L') = \Im(\Hom_{{\Cal O}_{X}}(L, L'\otimes _{
{\Cal O}_{X}}F_{p}{\Cal D}_{X}) \rightarrow  \Hom_{{\Bbb C}}(L, L')).
$$
We can verify that
$ F_{p} \Diff _{X}(L, L') $ coincides with the group of differential 
operators of order  $ \le  p $ in the sense of Grothendieck [6].
See [11, (1.20.2)].
But this filtration is not exhaustive in general.
(For example,
$ {\Cal D}_{X} \rightarrow  {\Cal O}_{X}  $ (defined by
$ P \mapsto  P1) $ has unbounded order if the right
$ {\Cal D}_{X} $-Module structure on
$ {\Cal D}_{X} $ is used.)

We denote by
$ M({\Cal O}_{X}, \Diff) $ the category such that the objects are
$ {\Cal O}_{X} $-Modules and the morphisms are given by
$  \Diff _{X}(L,L') $.
It is an additive category.
We have a functor
$$
\DR_{X}^{-1} : M({\Cal O}_{X}, \Diff) \rightarrow  M({\Cal D}_{X})^{r}
$$
such that
$ \DR_{X}^{-1}(L) = L\otimes _{{\Cal O}_{X}}{\Cal D}_{X} $,
using the injectivity of (1.2.2).

We can define the category of complexes
$ C({\Cal O}_{X}, \Diff) $, and also the category
$ K({\Cal O}_{X}, \Diff) $ whose morphisms are considered up to 
homotopy as in [12].
(An object of
$ C({\Cal O}_{X}, \Diff) $ or
$ K({\Cal O}_{X}, \Diff) $ is called a differential complex in this paper.)
Then
$ \DR_{X}^{-1} $ is extended to
$$
\DR_{X}^{-1} : C({\Cal O}_{X}, \Diff) \rightarrow  C({\Cal D}_{X})^{r},
\quad \DR_{X}^{-1} : K({\Cal O}_{X}, \Diff) \rightarrow  
K({\Cal D}_{X})^{r}.
$$

We say that
$ L^{\ssbull } \in  K({\Cal O}_{X}, \Diff) $ is
$ D $-acyclic if
$ \DR_{X}^{-1}(L^{\ssbull }) $ is acyclic, and a morphism
$ u : L^{\ssbull } \rightarrow  L^{\prime \ssbull } $ of
$ K({\Cal O}_{X}, \Diff) $ is a
$ D $-quasi-isomorphism if
$ \DR_{X}^{-1}(u) $ is a quasi-isomorphism.
Then by inverting
$ D $-quasi-isomorphisms, we get
$ D({\Cal O}_{X}, \Diff) $.
(We can verify that a
$ D $-quasi-isomorphism is a quasi-isomorphism by showing that 
(1.2.6) is a quasi-isomorphism.)
Similarly, we can define
$ C^{*}({\Cal O}_{X}, \Diff) $,  
$ K^{*}({\Cal O}_{X}, \Diff) $, 
$ D^{*}({\Cal O}_{X}, \Diff) $ for
$ * = +, -, b $ as in [12].

For a left
$ {\Cal D}_{X} $-Module
$ M $,
the de Rham complex
$ \DR_{X}(M) \in  C^{b}({\Cal O}_{X}, \Diff) $ is defined in the usual way 
(by identifying a left
$ {\Cal D}_{X} $-Module with an
$ {\Cal O}_{X} $-Module endowed with an integrable connection).
This gives a functor
$$
\DR_{X} : D^{*}({\Cal D}_{X})^{l} \rightarrow  D^{*}({\Cal O}_{X}, \Diff).
$$

For a right
$ {\Cal D}_{X} $-Modules
$ M $,
we define the de Rham complex
$ \DR_{X}(M) \in  C^{b}({\Cal O}_{X}, \Diff) $ so that the
$ i $-th component is
$ M\otimes _{{\Cal O}_{X}}\wedge ^{n-i}\Theta _{X} $ and the 
differential is given in the usual way.
(Taking a local coordinate system
$ (x_{1}, \dots , x_{n}) $,
it is isomorphic to the Koszul complex for the operators
$ \partial /\partial x_{1}, \dots , \partial /\partial x_{n} $ on
$ M $.)
This induces a functor
$$
\DR_{X} : D^{*}({\Cal D}_{X})^{r} \rightarrow  D^{*}({\Cal O}_{X}, \Diff),
$$
so that for
$ M^{\ssbull } \in  C({\Cal D}_{X})^{l} $ we have
$$
\DR_{X}(M^{\ssbull }) = \DR_{X}(\omega _{X}\otimes _{{\Cal O}_{X}}
M^{\ssbull })\quad  \text{in }D^{*}({\Cal O}_{X}, \Diff),
\tag 1.2.3
$$
using the transformation (1.1.1).

We have an equivalence of categories
$$
\DR_{X}^{-1} : D^{*}({\Cal O}_{X}, \Diff) \rightarrow  
D^{*}({\Cal D}_{X})^{r}
\tag 1.2.4
$$
for
$ * = +, -, b $ or empty, and a quasi-inverse is given by
$ \DR_{X}[\dim X] $.
This follows from the quasi-isomorphism
$$
\DR_{X}^{-1}\DR_{X}(M^{\ssbull })[\dim X] \rightarrow  M^{\ssbull }
\tag 1.2.5
$$
for
$ M^{\ssbull } \in  C({\Cal D}_{X})^{r} $.
See [11, (1.5.2)] and also Remark below.
Note that (1.2.5) implies a
$ D $-quasi-isomorphism
$$
\DR_{X}\DR_{X}^{-1}(L^{\ssbull })[\dim X] \rightarrow  L^{\ssbull }
\tag 1.2.6
$$
for
$ L^{\ssbull } \in  C({\Cal O}_{X}, \Diff) $ so that
$ \DR_{X}[\dim X] $ gives a quasi-inverse of
$ \DR_{X}^{-1} $.

\bigskip

\noindent
{\it Remark.} For the proof of (1.2.5) we consider the filtration on
$ \DR_{X}^{-1}\DR_{X}(M^{\ssbull }) $ induced by
$ F $ on
$ {\Cal D}_{X} $.
If we take a local coordinate system
$ (x_{1}, \dots , x_{n}) $,
then the direct sum of its graded pieces is the Koszul complex
$$
K(M\otimes _{{\Cal O}_{X}}\Gr^{F}{\Cal D}_{X}; \Gr^{F} \partial /\partial 
x_{1}, \dots , \Gr^{F} \partial /\partial x_{n})
$$
So (1.2.5) is clear.
(Here a differential operator of unbounded order is used.)
Similarly we can show that (1.2.6) is a quasi-isomorphism.

\bigskip
\noindent
{\bf 1.3.}
{\it Direct Images.} Let
$ f : X \rightarrow  Y $ be a morphism of smooth complex algebraic 
varieties.
For
$ M^{\ssbull } \in  D^{+}({\Cal D}_{X})^{r} $,
the direct image
$ f_{+}M^{\ssbull } \in  D^{+}({\Cal D}_{Y})^{r} $ is defined by
$$
f_{+}M^{\ssbull } = \bold{R}f_{*}(M^{\ssbull }{\otimes }_{{D}_{X}}
^{\bold{L}}{\Cal D}_{X\rightarrow Y}).
\tag 1.3.1
$$
Here
$ {\Cal D}_{X\rightarrow Y} = {\Cal O}_{X}\otimes _{f^{-1}
{\Cal O}_{Y}}f^{-1}{\Cal D}_{Y} $,
and it has a structure of left
$ {\Cal D}_{X} $-Module and right
$ f^{-1}{\Cal D}_{Y} $-Module.

For an
$ {\Cal O}_{X} $-Module
$ L $,
we define
$$
\DR_{Y}^{-1}(L) = \DR_{X}^{-1}(L)\otimes _{{\Cal D}_{X}}
{\Cal D}_{X\rightarrow Y} 
\,(= L\otimes _{f^{-1}{\Cal O}_{Y}}f^{-1}{\Cal D}_{Y}).
$$
Then  for
$ L^{\ssbull } \in  D^{+}({\Cal O}_{X}, \Diff) $, we have
$$
f_{+}\DR_{X}^{-1}(L^{\ssbull }) = \bold{R}f_{*}\DR_{Y}^{-1}(L^{\ssbull }) 
= \DR_{Y}^{-1}\bold{R}f_{*}(L^{\ssbull })\quad  \text{in }D^{+}(
{\Cal D}_{Y})^{r},
\tag 1.3.2
$$
where
$ \bold{R}f_{*}(L^{\ssbull }) $ is defined by taking a resolution of (a 
representative of)
$ L^{\ssbull } $ in
$ C^{+}({\Cal O}_{X}, \Diff) $ whose components are
$ f_{*} $-acyclic.
The isomorphisms of (1.3.2) mean that the differential complexes are 
stable by the sheaf-theoretic direct image, and
$ f_{+}\DR_{X}^{-1} = \DR_{Y}^{-1}\bold{R}f_{*} $.

For
$ M^{\ssbull } \in  D^{+}({\Cal D}_{X})^{l} $,
the direct image
$ f_{+}(M^{\ssbull }) $ is defined so that
$$
f_{+}(\omega _{X}\otimes _{{\Cal O}_{X}}M^{\ssbull }) = \omega 
_{Y}\otimes _{{\Cal O}_{Y}}f_{+}(M^{\ssbull })\quad  \text{in }
D^{+}({\Cal D}_{Y})^{r}.
\tag 1.3.3
$$
More explicitly, we have
$$
f_{+}M^{\ssbull } = \bold{R}f_{*}({\Cal D}_{Y\leftarrow X}{\otimes }
_{{D}_{X}}^{\bold{L}}M^{\ssbull }).
\tag 1.3.4
$$
Here
$ {\Cal D}_{Y\leftarrow X} $ is obtained by transforming the left and 
right
$ {\Cal D} $-Module structures on
$ {\Cal D}_{X\rightarrow Y} $ (by using the transformation of left and 
right
$ {\Cal D} $-Modules in (1.1.1)).

\proclaim{1.4.~Proposition}
Let  $ f : X \rightarrow  Y $ be a smooth morphism of smooth complex 
algebraic varieties with relative dimension
$ d $.
Then for
$ M^{\ssbull } \in  D^{+}({\Cal D}_{X})^{l} $,
we have a natural isomorphism of left
$ {\Cal D}_{Y} $-Modules :
$$
R^{i+d}f_{*}\DR_{X/Y}(M^{\ssbull }) = {\Cal H}^{i}(f_{+}M^{\ssbull }),
\tag 1.4.1
$$
where the left-hand side has the structure of left
$ {\Cal D}_{Y} $-Module by the Gauss-Manin connection [8].
\endproclaim

\demo\nofrills {Proof.\usualspace}
 Recall first the construction of Gauss-Manin connection in [8, \S 2].
Consider the Leray filtration
$ G $ on
$ {\Omega }_{X}^{i} $ defined by
$$
G^{p}{\Omega }_{X}^{i} = \Im(f^{*}{\Omega }_{Y}^{p}\otimes {\Omega }
_{X}^{i-p} \rightarrow  {\Omega }_{X}^{i}).
$$
Then
$ \Gr_{G}^{p}{\Omega }_{X}^{i} = f^{*}{\Omega }_{Y}^{p}\otimes _{
{\Cal O}_{X}}{\Omega }_{X/Y}^{i-p} $,
and
$ {\Omega }_{X/Y}^{i-p} $ is a locally free
$ {\Cal O}_{X} $-Module, because
$ f $ is smooth.
This induces a filtration
$ G $ on
$ \DR_{X}(M^{\ssbull }) $ such that
$$
\Gr_{G}^{p}\DR_{X}(M^{\ssbull }) = f^{*}{\Omega }_{Y}^{p}\otimes _{
{\Cal O}_{X}}\DR_{X/Y}(M^{\ssbull })[-p],
$$
where the
$ i $-th component of
$ \DR_{X/Y}(M^{j}) $ is
$ {\Omega }_{X/Y}^{j}\otimes _{{\Cal O}_{X}}M^{j} $.
So we get a spectral sequence
$$
{E}_{1}^{p,q} = {\Omega }_{Y}^{p}\otimes _{
{\Cal O}_{Y}}R^{q}f_{*}\DR_{X/Y}(M^{\ssbull }) \Rightarrow  
R^{p+q}f_{*}\DR_{X}(M^{\ssbull }).
\tag 1.4.2
$$
This is defined in the category of sheaves on
$ Y $.
However, the
$ E_{1} $-complex is defined in
$ C^{*}({\Cal O}_{Y}, \Diff) $, and is isomorphic to
$ \DR_{Y}(R^{q}f_{*}\DR_{X/Y}(M^{\ssbull })) $.
See loc.~cit.

The Gauss-Manin connection on
$ R^{q}f_{*}\DR_{X/Y}(M^{\ssbull }) $ is defined as
$ d_{1} : {E}_{1}^{0,q} \rightarrow  {E}_{1}^{1,q} $.
More generally,
$ d_{1} : {E}_{1}^{p,q} \rightarrow  {E}_{1}^{p+1,q} $ is given by the 
connecting morphism associated with the short exact sequence
$$
0 \rightarrow  \Gr_{G}^{p+1}\DR_{X}(M^{\ssbull }) \rightarrow  
G^{p}\DR_{X}(M^{\ssbull })/G^{p+2}\DR_{X}(M^{\ssbull }) \rightarrow  
\Gr_{G}^{p}\DR_{X}(M^{\ssbull }) \rightarrow  0.
\tag 1.4.3
$$

The filtration
$ G $ induces also a filtration
$ G $ on
$ \DR_{Y}^{-1}\DR_{X}(M^{\ssbull }) $ in the notation of (1.3.2), and
$ \DR_{Y}^{-1} $ commutes with
$ \Gr^{G} $,
because
$ {\Cal D}_{Y} $ is flat over
$ {\Cal O}_{Y} $.
So we get a spectral sequence defined in the category of right
$ {\Cal D}_{Y} $-Modules :
$$
{E}_{1}^{p,q} = \DR_{Y}^{-1}({\Omega }_{Y}^{p}\otimes _{
{\Cal O}_{Y}}R^{q}f_{*}\DR_{X/Y}(M^{\ssbull })) \Rightarrow  
{\Cal H}^{p+q}(\DR_{Y}^{-1}\bold{R}f_{*}\DR_{X}(M^{\ssbull })).
\tag 1.4.4
$$
Applying (1.2.1) to the
$ E_{1} $-complexes of (1.4.2) and (1.4.4) (and using (1.4.3)), we see 
that the
$ E_{1} $-complex of (1.4.4) is isomorphic to
$ \DR_{Y}^{-1}\DR_{Y}(R^{q}f_{*}\DR_{X/Y}(M^{\ssbull })) $.
The latter is further quasi-isomorphic to
$$
\omega _{Y}\otimes _{{\Cal O}_{Y}}(R^{q}f_{*}\DR_{X/Y}(M^{\ssbull }))
[-m]
$$
by the quasi-isomorphism (1.2.5) combined with (1.2.3), where
$ m = \dim Y $.
This implies that the spectral sequence (1.4.4) degenerates at
$ E_{2} $,
and we get the isomorphism of right
$ {\Cal D}_{Y} $-Modules :
$$
\omega _{Y}\otimes _{{\Cal O}_{Y}}(R^{q}f_{*}\DR_{X/Y}(M^{\ssbull })) = 
\omega _{Y}\otimes _{{\Cal O}_{Y}}{\Cal H}^{q-d}f_{+}(M^{\ssbull }),
$$
because
$$
{\Cal H}^{i}(\DR_{Y}^{-1}\bold{R}f_{*}\DR_{X}(M^{\ssbull })) = \omega 
_{Y}\otimes _{{\Cal O}_{Y}}{\Cal H}^{i-n}f_{+}(M^{\ssbull })
$$
by (1.2.5), (1.2.3), (1.3.2) and (1.3.3), where
$ n = \dim X $.
So the assertion follows.
\enddemo

\bigskip
\noindent
{\bf 1.5.}
{\it Pull-backs and local cohomology.} Let
$ f : X \rightarrow  Y $ be a morphism of smooth complex algebraic 
varieties.
For
$ M \in  M({\Cal D}_{Y})^{l} $,
the usual pull-back
$ f^{*}M \,(:= {\Cal O}_{X}\otimes _{f^{-1}{\Cal O}_{Y}} f^{-1}M) $ has 
naturally a left
$ {\Cal D}_{X} $-Module structure by
$$
f^{*}M = {\Cal D}_{X\rightarrow Y}\otimes _{f^{-1}{\Cal D}_{Y}} f^{-1}M,
\tag 1.5.1
$$
where
$ {\Cal D}_{X\rightarrow Y} $ is as in (1.3).
So, using a flat resolution, we get the derived functor
$$
\bold{L}f^{*} : D^{*}({\Cal D}_{Y})^{l} \rightarrow  D^{*}({\Cal D}_{X})^{l}
\tag 1.5.2
$$
for
$ * = +, -, b $ or empty.
This is compatible with
$ \bold{L}f^{*} $ for
$ {\Cal O} $-Modules, because a flat
$ {\Cal D}_{X} $-Module is flat over
$ {\Cal O}_{X} $.
We will use the notation
$ f^{*}M^{\ssbull } $ for
$ \bold{L}f^{*}M^{\ssbull } $ if
$ f $ is flat, or more generally, if
$ f $ is cohomologically noncharacteristic for
$ M^{\ssbull } $ in the sense that the higher torsion groups of
$ {\Cal O}_{X} $ and
$ f^{-1}M^{i} $ over
$ f^{-1}{\Cal O}_{Y} $ vanish.

If we identify a left
$ {\Cal D}_{Y} $-Module with an
$ {\Cal O}_{Y} $-Module (not necessarily finitely generated) endowed 
with an integrable connection
$ \nabla  $,
the pull-back of
$ {\Cal D}_{Y} $-Modules corresponds to the pull-back of connections 
(i.e.,
$ \nabla  : M \rightarrow  {\Omega }_{Y}^{1}\otimes _{{\Cal O}_{Y}}M $ 
induces
$ f^{*}\nabla  : f^{*}M \rightarrow  {\Omega }_{X}^{1}\otimes _{
{\Cal O}_{X}}f^{*}M) $.

We define also
$$
f^{!}M^{\ssbull } = \bold{L}f^{*}M^{\ssbull }[\dim X - \dim Y].
\tag 1.5.3
$$
This
$ f^{!} $ is compatible with
$ f^{!} $ for
$ {\Cal O} $-Modules up to tensor of
$ \omega _{X/Y} $.
More precisely, if we define
$ f^{!} $ for right
$ {\Cal D} $-Modules using the transformation of left and right
$ {\Cal D} $-Modules as in (1.1.1), then this is compatible with
$ f^{!} $ for
$ {\Cal O} $-Modules.

Let
$ i : Z \rightarrow  X $ be a closed embedding of smooth complex 
algebraic varieties.
Then we have
$$
i_{+}i^{!}M^{\ssbull } = \bold{R}\Gamma _{Z}M^{\ssbull }
\tag 1.5.4
$$
for
$ M^{\ssbull } \in  {D}_{\qcoh}^{+}({\Cal D}_{X})^{l} $,
where
$ \Gamma _{Z}M $ for a sheaf
$ M $ is the subsheaf consisting of local sections supported in
$ Z $,  and
$  \bold{R}\Gamma _{Z}M^{\ssbull } $ is defined by taking a flasque resolution.
For the proof of (1.5.4), we consider a functor
$ \bold{R}\Gamma _{[Z]} $ defined by
$$
\bold{R}\Gamma _{[Z]}M^{\ssbull } = 
\hbox{{\rm ind lim}}_{j} {\cHom}
_{{\Cal O}_{X}}({\Cal O}_{X}/{I}_{Z}^{j}, M^{\prime\ssbull }), 
$$
where
$ I_{Z} $ is the ideal sheaf of
$ Z $, and
$ M^{\ssbull } \to M^{\prime\ssbull } $ is a quasi-isomorphism such that
$ M^{\prime\ssbull } $ is a complex of injective
$ {\Cal D}_{X} $-Modules.
See for example [2, VI, 7.9].
Then for 
$ M^{\ssbull } \in  {D}_{\qcoh}^{+}({\Cal D}_{X}) $,  we have a natural
isomorphism
$$
\bold{R}\Gamma _{[Z]}M^{\ssbull } \simto
\bold{R}\Gamma _{Z}M^{\ssbull }
\tag 1.5.5
$$
in the derived category.
Furthermore
$ i^{!}M^{\ssbull } $ for right
$ {\Cal D} $-Modules is given by
$$
i^{!}M^{\ssbull } = \bold{R}{\cHom}_{{\Cal O}_{X}}
({\Cal O}_{X}/I_{Z}, M^{\ssbull }).
$$
Indeed,
$ \Hom_{{\Cal O}_{X}}({\Cal O}_{X}/I_{Z}, M) $ has naturally a structure of right
$ {\Cal D}_{Z} $-Module if
$ M $ is a right
$ {\Cal D}_{X} $-Module.
So we get a canonical morphism in
$ {D}_{\qcoh}^{+}({\Cal D}_{X}) $
$$
i_{+}i^{!}M^{\ssbull } \rightarrow  \bold{R}\Gamma _{[Z]}M^{\ssbull }
\tag 1.5.6
$$
and it is enough to show that this is an isomorphism.
Then we may assume that 
$ M^{\ssbull } $ is a quasi-coherent
$ {\Cal D}_{X} $-Module
$ M $,
and the assertion follows by taking a resolution of
$ M $ such that the components are quasi-coherent
$ {\Cal D}_{X} $-Modules, and are injective objects in the category of 
quasi-coherent
$ {\Cal O}_{X} $-Modules.
See [2, VI, 7.13].

\proclaim{1.6.~Proposition}
Let  $ f : X \rightarrow  Y $ be a morphism of smooth complex algebraic 
varieties, and
$ Z' $ a closed subvariety of
$ Y $.
Let
$ Z = f^{-1}(Z') $,  and consider a cartesian diagram
$$
\CD
Z @>i>>   X @<j<< X\backslash Z
\\
@VVf'V    @VVfV @ VVf''V
\\
Z'   @>i'>>   Y @<j'<< Y\backslash Z'
\endCD
$$
where the horizontal morphisms are inclusions 
and the vertical ones are induced by 
$ f $.
Then for
$ M^{\ssbull } \in  {D}_{\qcoh}^{+}({\Cal D}_{X})^{l} $,
we have a canonical isomorphism in
$ {D}_{\qcoh}^{+}({\Cal D}_{Y})^{l} $
$$
f_{+}\bold{R}\Gamma _{Z}M^{\ssbull } = \bold{R}\Gamma 
_{Z'}f_{+}M^{\ssbull }.
\tag 1.6.1
$$
If
$ Z $ and
$ Z' $ are smooth, we have as in [2, VI, 8.4] :
$$
f'_{+}i^{!}M^{\ssbull } = i^{\prime !}f_{+}M^{\ssbull }.
\tag 1.6.2
$$
\endproclaim

\demo\nofrills {Proof.\usualspace}
We have a distinguished triangle in
$ {D}_{\qcoh}^{+}({\Cal D}_{X})^{l} $
$$
\rightarrow  \bold{R}\Gamma _{Z}M^{\ssbull } \rightarrow  M^{\ssbull } 
\rightarrow  j_{+}j^{*}M^{\ssbull } \rightarrow ,
$$
because
$  j_{+}j^{*}M^{\ssbull } = \bold{R}j_{*}j^{*}M^{\ssbull } $ by definition.
Applying the functorial morphism
$ \bold{R}\Gamma _{Z'}f_{+} \rightarrow  f_{+}  $ to this triangle, 
(1.6.1) is reduced to the following two assertions for
$ M^{\prime \ssbull } \in  {D}_{\qcoh}^{+}({\Cal D}_{Y})^{l} $:
$$
\align
&\bold{R}\Gamma _{Z'}M^{\prime \ssbull } \simto M^{\prime \ssbull }\,\,\, 
\text{if } \supp  M' \subset  Z',
\\
&\bold{R}\Gamma _{Z'}j'_{+}j^{\prime *}M^{\prime \ssbull } = 0,
\endalign
$$
where
$ M^{\prime \ssbull } $ is respectively
$ f_{+}\bold{R}\Gamma _{Z}M^{\ssbull } $ and
$ f_{+}M^{\ssbull } $ (because
$ j'_{+}j^{\prime *}f_{+}M^{\ssbull } = j'_{+}f''_{+}j^{*}M^{\ssbull } = 
f_{+}j_{+}j^{*}M^{\ssbull }) $.
The first assertion is clear.
For the second, it is enough to apply the functor
$ \bold{R}\Gamma _{Z'} $ to the distinguished triangle
$$
\rightarrow  \bold{R}\Gamma _{Z'}M^{\prime \ssbull } \rightarrow  M^{\prime \ssbull }
 \rightarrow  j'_{+}j^{\prime *}M^{\prime \ssbull } \rightarrow .
$$

If
$ Z $ and
$ Z' $ are smooth, (1.6.2) follows from (1.6.1) by the Kashiwara 
equivalence (see for example [2, VI, 7.13]).
\enddemo

\noindent
{\it Remark.} The above argument is more precise than in [2, VI, 8.4] 
about the problem of the ambiguity of mapping cone.

\bigskip
\noindent
{\bf 1.7.}
{\it Twists.} Let
$ X $ be a smooth complex algebraic variety, and
$ F \in  \Gamma (X, {\Cal O}_{X}) $.
Then
$ {\Cal O}_{X}e^{F} $ is defined to be the quotient of
$ {\Cal D}_{X} $ by the left ideal generated by
$ \xi  - \xi (F) $ for vector fields
$ \xi  $ (where
$ \xi (F) \in  {\Cal O}_{X}) $.
As an
$ {\Cal O}_{X} $-Module,
$ {\Cal O}_{X}e^{F} $ is isomorphic to
$ {\Cal O}_{X} $ (generated by
$ 1 \in  {\Cal D}_{X}) $.
For a left
$ {\Cal D}_{X} $-Module
$ M $,
we define the twist
$ M_{F} $ (or
$ M\otimes e^{F}) $ by
$$
M_{F} = M\otimes _{{\Cal O}_{X}}{\Cal O}_{X}e^{F}.
$$
See also [9].
In terms of integrable connection,
$ M_{F} $ corresponds to the twist of connection
$ \nabla _{F} $ defined by
$ \nabla _{F}(m) = \nabla m + dF\otimes m $.
So we have
$$
\DR_{X}(M_{F}) = \DR_{X}(M, \nabla  + dF\wedge ).
\tag 1.7.1
$$

Let
$ g : Z \rightarrow  X $ be a morphism of smooth complex algebraic 
varieties.
Then we have
$$
g^{*}(M_{F}) = (g^{*}M)_{g^{*}F}.
\tag 1.7.2
$$

\bigskip\bigskip\centerline{{\bf 2.~Dwork Cohomology}}

\bigskip
\noindent
{\bf 2.1.}
Let
$ X $ be a smooth complex algebraic variety, and
$ \pi  : V \rightarrow  X $ an algebraic vector bundle of rank
$ r $.
Let
$ \pi ^{\vee } : V^{\vee } \rightarrow  X $ be the dual vector bundle of
$ V $.
We denote by
$ {\Cal O}_{X}(V^{\vee }) $ the sheaf of sections of
$ V^{\vee } $,
and by
$ \Sym_{{\Cal O}_{X}}{\Cal O}_{X}(V^{\vee }) $ the symmetric algebra 
over
$ {\Cal O}_{X} $.
Then
$$
V =  \Spec_{X}(\Sym_{{\Cal O}_{X}}{\Cal O}_{X}(V^{\vee })).
$$
So a section
$ s $ of
$ V^{\vee } $ is identified with a function
$ F $ on
$ V $ such that the restriction to each fiber of
$ V $ is a linear function.

We consider a cartesian diagram
$$
\CD
V @>\tilde{s}>>   V\times _{X}V^{\vee }   @>\tilde{\pi }^{\vee}>>   V
\\
@VV{\pi}V       @VV{\tilde{\pi }}V       @VV{\pi}V 
\\
X   @>s>>   V^{\vee }  @>{\pi^{\vee}}>>  X
\endCD
\tag 2.1.1
$$
Let
$ \langle v, v^{\vee } \rangle \in  {\Bbb C} $ denote the canonical 
pairing of
$ v \in  V, v^{\vee } \in  V^{\vee } $ such that
$ \pi (v) = \pi ^{\vee }(v^{\vee }) $.
This gives a function
$ \Phi  \in  \Gamma (V\times _{X}V^{\vee },{\Cal O}) $ such that
$ \Phi (v, v^{\vee }) = \langle v, v^{\vee } \rangle $ for
$ (v, v^{\vee }) \in  V\times _{X}V^{\vee } $.

Let
$ (y_{1}, \dots , y_{r}), ({y}_{1}^{\vee }, \dots , {y}_{r}^{\vee }) $ be 
locally defined dual coordinate systems of the vector bundles
$ V $ and
$ V^{\vee } $.
This means that we have dual local bases
$ (e_{1}, \dots , e_{r}), ({e}_{1}^{\vee }, \dots , {e}_{r}^{\vee }) $ of
$ {\Cal O}_{X}(V), {\Cal O}_{X}(V^{\vee }) $ such that
$ \sum _{i} y_{i} (v) e_{i}(\pi (v))  = v $ for
$ v \in  V $ (and similarly for
$ V^{\vee }) $.
Then we have
$ \Phi (v, v^{\vee }) = \sum _{i} y_{i} (v) {y}_{i}^{\vee }(v^{\vee }) $ 
for
$ (v, v^{\vee }) \in  V\times _{X}V^{\vee } $,
and
$$
F = \tilde{s}^{*}\Phi .
\tag 2.1.2
$$
Indeed,
$ F = \sum _{i }y_{i}f_{i} $ if
$ s $ is expressed by
$ (f_{1}, \dots , f_{r}) $ using the coordinates
$ ({y}_{1}^{\vee }, \dots , {y}_{r}^{\vee }) $.

\bigskip
\noindent
{\bf 2.2.}
{\it First proof of} (0.2).
Let
$ \tilde{\tilde{\pi }} : V\times _{X}V^{\vee } \rightarrow  X $ denote the 
natural projection, and
$ i^{\vee } : X \rightarrow  V^{\vee } $ the zero section.
We prove the assertion by showing the following isomorphisms :
$$
\bold{R}\Gamma _{Y}M^{\ssbull } = s^{!}{i}_{+}^{\vee }M^{\ssbull } = 
s^{!}\tilde{\pi }_{+}((\tilde{\tilde{\pi }}^{*}M^{\ssbull })_{\Phi }) = \pi 
_{+}\tilde{s}^{!}((\tilde{\tilde{\pi }}^{*}M^{\ssbull })_{\Phi }) = \pi 
_{+}((\pi ^{*}M^{\ssbull })_{F})[-r].
\tag 2.2.1
$$

The first isomorphism follows from
$$
s^{!}{i}_{+}^{\vee }M^{\ssbull } = {\pi }_{+}^{\vee }s_{+}s^{!}{i}_{+}^{\vee }
M^{\ssbull } = {\pi }_{+}^{\vee }\bold{R}\Gamma _{s(X)}{i}_{+}^{\vee }
M^{\ssbull } = {\pi }_{+}^{\vee }{i}_{+}^{\vee }\bold{R}\Gamma 
_{Y}M^{\ssbull } = \bold{R}\Gamma _{Y}M^{\ssbull },
$$
using (1.5.4) and (1.6).

As for the last isomorphism of (2.2.1), note that
$ \tilde{s} $ is noncharacteristic for
$ (\tilde{\tilde{\pi }}^{*}M^{\ssbull })_{\Phi } $ (or equivalently for
$ \tilde{\tilde{\pi }}^{*}M^{\ssbull } = (\pi ^{\vee })^{*}\pi ^{*}M^{\ssbull }) $,
because it is a section of a smooth morphism
$ \pi ^{\vee } $.
So we will write
$ \tilde{s}^{*} $ for
$ \bold{L}\tilde{s}^{*} $.
See (1.5).
Then
$$
(\pi ^{*}M^{\ssbull })_{F} = \tilde{s}^{*}((\tilde{\tilde{\pi }}
^{*}M^{\ssbull })_{\Phi }) = 
\tilde{s}^{!}((\tilde{\tilde{\pi }}^{*}M^{\ssbull })_{\Phi }[r]
$$
by (1.7.2) and by definition (i.e.
 $ \tilde{s}^{!} = \tilde{s}^{*}[-r]) $.

The third isomorphism of (2.2.1) follows from the base change
$ s^{!}\tilde{\pi }_{+} = \pi _{+}\tilde{s}^{!} $ (see (1.6)).
Then it remains to show the second isomorphism of (2.2.1), and the 
assertion is reduced to the following.

\proclaim{2.3.~Lemma}
 $ {i}_{+}^{\vee }M^{\ssbull } = \tilde{\pi }_{+}((\tilde{\tilde{\pi }}
^{*}M^{\ssbull })_{\Phi }) $.
\endproclaim

\demo\nofrills {Proof.\usualspace}
 We may assume
$ M^{\ssbull } $ is a quasi-coherent left
$ {\Cal D}_{X} $-Module
$ M $,
because the assertion in this case implies
$ {\Cal H}^{i}\tilde{\pi }_{*}((\tilde{\tilde{\pi }}^{*}M^{j})_{\Phi }) = 0 $ 
for
$ i \ne  0 $.
Note that the right-hand side is the Fourier transform of
$ \pi ^{*}M $,
and the assertion is more or less well-known if the vector bundle
$ V $ is trivial.
Indeed, the right-hand side in this case is the (shifted) Koszul complex
$$
K(M[y_{1}, \dots , y_{r},{y}_{1}^{\vee }, \dots , {y}_{r}^{\vee }] ; 
\partial /\partial y_{i} + {y}_{i}^{\vee } \,(1\le  i \le  r))[r]
\tag 2.3.1
$$
by definition of twist (1.7), where
$ (y_{1}, \dots , y_{r}), ({y}_{1}^{\vee }, \dots , {y}_{r}^{\vee }) $ are as 
in (2.1).

In general it is enough to show that the local isomorphisms give a 
globally well-defined isomorphism.
So the assertion is reduced to
$$
M = \tilde{\tilde{\pi }}_{+}((\tilde{\tilde{\pi }}^{*}M)_{\Phi }).
\tag 2.3.2
$$

Let
$ \tilde{V} = V\times _{X}V^{\vee } $.
By definition,
$ {\Cal H}^{0}\tilde{\tilde{\pi }}_{+}((\tilde{\tilde{\pi }}^{*}M)_{\Phi }) 
$ is a quotient of
$ \tilde{\tilde{\pi }}_{*}({\Omega }_{\tilde{V}/X}^{2r}\otimes _{
{\Cal O}_{V}}\tilde{\tilde{\pi }}^{*}M) $ which is isomorphic to
$ \tilde{\tilde{\pi }}_{*}{\Omega }_{\tilde{V}/X}^{2r}\otimes _{
{\Cal O}_{X}}M $ because
$ M $ is quasi-coherent.
Then
$ {\Omega }_{\tilde{V}/X}^{2r} $ is globally trivialized by the relative 
form
$$
dy_{1}\wedge  \dots \wedge dy_{r}\wedge d{y}_{1}^{\vee }\wedge  
\dots \wedge d{y}_{r}^{\vee },
$$
because this relative form is independent of the choice of the locally 
defined dual coordinate systems
$ (y_{1}, \dots , y_{r}), ({y}_{1}^{\vee }, \dots , {y}_{r}^{\vee }) $ of the 
vector bundles
$ V, V^{\vee } $.
So
$ {\Cal H}^{0}\tilde{\tilde{\pi }}_{+}((\tilde{\tilde{\pi }}^{*}M)_{\Phi }) 
$ is globally a quotient of
$ \tilde{\tilde{\pi }}_{*}{\Cal O}_{\tilde{V}}\otimes _{{\Cal O}_{X}}M $.
Since
$ \tilde{\tilde{\pi }}_{+}((\tilde{\tilde{\pi }}^{*}M)_{\Phi }) $ is locally 
the (shifted) Koszul complex
$$
K(M[y_{1}, \dots , y_{r,}{y}_{1}^{\vee }, \dots , {y}_{r}^{\vee }] ; 
\partial /\partial y_{i} + {y}_{i}^{\vee }, 
\partial /\partial {y}_{i}^{\vee } + y_{i} \,(1\le  i \le  r))[2r],
\tag 2.3.3
$$
we see that
$ M =1\otimes M \subset  \tilde{\tilde{\pi }}_{*}
{\Cal O}_{\tilde{V}}\otimes _{{\Cal O}_{X}}M $ is isomorphic to
$ {\Cal H}^{0}\tilde{\tilde{\pi }}_{+}((\tilde{\tilde{\pi }}^{*}M)_{\Phi }) 
$,
and this isomorphism is compatible with the action of
$ {\Cal D}_{X} $.
So the assertion follows.
\enddemo

\bigskip
\noindent
{\bf 2.4.}
{\it Second proof of} (0.2).
We first construct a canonical morphism in
$ C^{b}({\Cal O}_{X},  \Diff) $ :
$$
\pi _{*}(\DR_{V}(\pi ^{*}M^{\ssbull }), \nabla  + dF\wedge ) \rightarrow  
\DR_{X}(M^{\ssbull }).
\tag 2.4.1
$$
Note that the left-hand side is
$ \pi _{*}(\DR_{V}((\pi ^{*}M^{\ssbull })_{F}) $ by (1.7.1).
Since the
$ M^{j} $ are quasi-coherent, we have the projection formula
$$
\pi _{*}({\Omega }_{V}^{i}\otimes _{{\Cal O}_{V}}\pi ^{*}M^{j}) = (\pi 
_{*}{\Omega }_{V}^{i})\otimes _{{\Cal O}_{X}}M^{j}
$$
The natural
$ {\Bbb C}^{*} $-action on the vector bundle
$ V $ gives a decomposition
$$
\pi _{*}{\Omega }_{V}^{i} = \bigoplus _{k\in {\Bbb N}} (\pi _{*}{\Omega }
_{V}^{i})^{\langle k \rangle},
$$
such that
$ \lambda ^{*}\omega  = \lambda ^{k}\omega  $ for
$ \lambda  \in  {\Bbb C}^{*}, \omega  \in  (\pi _{*}{\Omega }
_{V}^{i})^{\langle k \rangle} $.
(This is easily verified by taking a local trivialization of
$ V $.)
So we get a decomposition
$$
\pi _{*}({\Omega }_{V}^{i}\otimes _{{\Cal O}_{V}}\pi ^{*}M^{j}) = 
\bigoplus _{k\in {\Bbb N}} (\pi _{*}{\Omega }_{V}^{i})^{\langle k 
\rangle}\otimes _{{\Cal O}_{X}}M^{j}.
$$
The differential
$ \nabla  $ preserves the decomposition, and
$ dF\wedge  $ preserves it up to a shift of degree by one because the 
restriction of
$ F $ on each fiber is linear.
Since
$ (\pi _{*}{\Omega }_{V}^{i})^{\langle 0 \rangle} = {\Omega }_{X}^{i} $,
 $ \DR_{X}(M^{\ssbull }) $ is a quotient complex of
$ \pi _{*}(\DR_{V}(\pi ^{*}M^{\ssbull }), \nabla  + dF\wedge ) $,
and we get (2.4.1).

Let
$ j : X \backslash  Y \rightarrow  X $ denote the inclusion morphism.
We have a distinguished triangle
$$
\rightarrow  \bold{R}\Gamma _{Y}M^{\ssbull } \rightarrow  M^{\ssbull } 
\rightarrow  \bold{R}j_{*}j^{*}M^{\ssbull } \rightarrow ,
$$
where
$ \bold{R}\Gamma _{Y}M^{\ssbull }, \bold{R}j_{*}j^{*}M^{\ssbull } $ are 
defined by using a flasque quasi-coherent resolution.
Applying (2.4.1) to this triangle, it is enough to show the following two 
assertions for the proof of (0.2) :
$$
\align
&\text{the source of (2.4.1) for 
$ \bold{R}j_{*}j^{*}M^{\ssbull } $ is  $ D $-acyclic},
\tag 2.4.2
\\
&\text{(2.4.1) is a $ D $-quasi-isomorphism if
$ \Supp M^{\ssbull } \subset  Y $.}
\tag 2.4.3
\endalign
$$
So we may assume that
$ M^{\ssbull } $ is a quasi-coherent left
$ {\Cal D}_{X} $-Module
$ M $.

Since the assertions are local, we may assume
$ V = X \times  {\Bbb A}^{r} $,
and
$ s = (f_{1}, \dots , f_{r}) $.
Taking
$ \DR_{X}^{-1} $ of (2.4.1) and then transforming right
$ {\Cal D} $-Modules to left
$ {\Cal D} $-Modules, (2.4.1) becomes a morphism of complex of left
$ {\Cal D}_{X} $-Modules :
$$
\pi _{*}(\DR_{V/X}(\pi ^{*}M), \nabla  + dF\wedge ) \rightarrow  M.
\tag 2.4.4
$$
So it is enough to show (2.4.2--3) with (2.4.1) replaced by (2.4.4), and
$ D $-acyclic and
$ D $-quasi-isomorphism by acyclic and quasi-isomorphism 
respectively.

Since
$ F = \sum _{i} y_{i}f_{i} $ with the notation of (2.1.2),  
the left-hand side of (2.4.4) is the Koszul complex
$$
K(M[y_{1}, \dots , y_{r}]; \partial /\partial y_{i} + f_{i} 
\,(1\le  i \le  r)).
\tag 2.4.5
$$
Let
$ U_{i} = \{f_{i} \ne  0\} \subset  X $,
and
$ U_{I} = \cap _{i\in I} U_{i} $ with the natural inclusion
$ j_{I} : U_{I} \rightarrow  X $ for
$ I \subset  \{1, \dots , r\} $.
Then we have a Cech complex with
$ i $-th component defined by
$$
\bigoplus _{|I|=j+1}(j_{I})_{*}{j}_{I}^{{}^*}M\,\,\, \text{for }i \ge  0,\,\,\, 
\text{and }0\,\,\, \text{otherwise} .
$$
Since this complex represents
$ \bold{R}j_{*}j^{*}M $,
we may replace
$ \bold{R}j_{*}j^{*}M^{\ssbull } $ in (2.4.2) by
$ (j_{I})_{*}{j}_{I}^{{}^*}M $ for
$ I \subset  \{1, \dots , r\} $.
By definition, the action of
$ f_{i} $ on
$ (j_{I})_{*}{j}_{I}^{{}^*}M $ is bijective for
$ i \in  I $.
This implies the bijectivity of the action of
$ \partial /\partial y_{i} + f_{i} $ on
$ (j_{I})_{*}{j}_{I}^{{}^*}M[y_{1}, \dots , y_{r}] $ using the increasing 
filtration on
$ {\Bbb C}[y_{1}, \dots , y_{r}] $ by degree.
So we get (2.4.2), because the Koszul complex (2.4.5) is the single 
complex associated with the
$ r $-ple complex defined by the morphisms
$ \partial /\partial y_{i} + f_{i} \,(1\le  i \le  r)  $.

For (2.4.3) we take an increasing exhaustive filtration
$ G $ of
$ M $ such that
$ G_{-1}M = 0 $ and
$ \Gr_{i}^{G}M $ are annihilated by the reduced ideal sheaf of
$ Y $.
Then
$ \Gr^{G} $ of (2.4.5) is the Koszul complex
$$
K(\Gr^{G}M[y_{1}, \dots , y_{r}]; \partial /\partial y_{i} \,(1\le  i \le  r)),
$$
and
$ \Gr^{G} $ of (2.4.4) is a quasi-isomorphism.
So we get (2.4.3).
This completes the second proof of (0.2).

\bigskip
\noindent
{\bf 2.5.}
{\it Remarks.} (i) The cohomology sheaf on the right-hand side of (0.2) 
is called the generalized Dwork cohomology sheaf.
If
$ p $ is smooth of relative dimension
$ n $,
we have by (1.4) an isomorphism of left
$ {\Cal D}_{Y} $-Modules :
$$
{\Cal H}^{i}(p\scirc \pi )_{+}((\pi ^{*}M^{\ssbull })_{F}) = 
R^{i+n+r}(p\scirc \pi )_{*}\DR_{V/S}((\pi ^{*}M^{\ssbull })_{F}).
\tag 2.5.1
$$
If
$ M^{\ssbull } $ is a coherent
$ {\Cal O}_{X} $-Module with an integrable connection
$ ({\Cal E},\nabla ) $, then
$ (\pi ^{*}M^{\ssbull })_{F} $ is the pull-back
$ \pi ^{*}{\Cal E} $ with the twisted connection
$ \nabla  + dF\wedge  $.
If furthermore
$ Y \rightarrow  S $ is smooth with
$ \codim _{X}Y = d $,
then
$ \bold{R}\Gamma _{Y}M^{\ssbull }[d] = i_{+}i^{*}({\Cal E},\nabla ) $ 
where
$ i : Y \rightarrow  X $ denotes the inclusion morphism.
So we get
$$
{\Cal H}^{i}(p_{+}\bold{R}\Gamma _{Y}M^{\ssbull }[d]) = 
R^{i+n-d}(p\scirc i)_{*}\DR_{Y/S}(i^{*}({\Cal E},\nabla )).
\tag 2.5.2
$$
Then (0.3) follows from (0.2) combined with (2.5.1) and (2.5.2).

(ii) In the case
$ M^{\ssbull } = {\Cal O}_{X} $ and
$ Y $ is connected, we can show that the isomorphisms obtained by the 
two proofs of (0.2) coincide up to a nonzero constant multiple using the 
Riemann-Hilbert correspondence together with the duality.
Indeed, we may assume that
$ p $ is the identity,
and we have an isomorphism
$$
\Hom(\bold{R}\Gamma _{Y}{\Cal O}_{X}, \bold{R}\Gamma _{Y}
{\Cal O}_{X}) = \Hom({\Bbb C}_{Y^{\an}}, {\Bbb C}_{Y^{\an}}) = {\Bbb C},
\tag 2.5.3
$$
where the first
$ \Hom $ is taken in the derived category of left
$ {\Cal D}_{X} $-Modules, and the second
$ \Hom $ in the derived category of
$ {\Bbb C}_{X^{\an}} $-Modules.

\bigskip\bigskip\centerline{{\bf 3.~Projective Hypersurface Case}}

\bigskip
\noindent
{\bf 3.1.}
With the notation of (2.1), assume
$ X = {\Bbb P}_{S}^{n} $ with
$ S $ a smooth connected variety, and
$ p : X \rightarrow  S $ is the natural projection.
(For example,
$ S = {\Bbb A}^{r} $.
Here we do not have to restrict to an open subvariety of
$ {\Bbb A}^{r} $,
because
$ Y $ is not assumed to be smooth over
$ S $.)

For an integer
$ a $,
we define a line bundle
$ \pi  : V^{(a)} \rightarrow  X $ by
$$
V^{(a)} =  \Spec _{X}(\bigoplus _{k\in {\Bbb N}} {\Cal O}_{X}(ak))
\tag 3.1.1
$$
so that
$ {\Cal O}_{X}(-a) = {\Cal O}_{X}(V^{(a)})  $ (the sheaf of sections of
$ V^{(a)}) $.
If
$ a = 1 $,
 $ V^{(1)} $ is the tautological line bundle, and is the blow-up of
$ {\Bbb A}_{S}^{n+1} $ at the origin.

We assume
$ V = V^{(m)} $.
Then a nonzero section
$ s $ of
$ {\Cal O}_{X}(m) $ is identified with a homogeneous polynomial of 
degree
$ m $ with coefficients in
$ {\Cal O}_{S} $,
which is denoted also by
$ F $.
(This is compatible with the notation
$ F \in  \Gamma (V^{(m)}, {\Cal O}_{V}(m)) $ in (2.1) by (3.2).)
For
$ M^{\ssbull } = {\Cal O}_{X} $,
we have the Dwork cohomology sheaf
$ R^{i}(p\scirc \pi )_{*}({\Omega }_{{V}^{(m)}/S}^{\ssbull }, d + 
dF\wedge ) $ as in the introduction.
It has a structure of left
$ {\Cal D}_{S} $-Module twisted by
$ F $ (which is compatible with (2.5.1)).
By (2.4) (combined with (2.5.1)), we get a canonical isomorphism of left
$ {\Cal D}_{S} $-Modules
$$
{\Cal H}^{i}p_{+}\bold{R}\Gamma _{Y}{\Cal O}_{X} = R^{i+n}(p\scirc \pi )
_{*}({\Omega }_{{V}^{(m)}/S}^{\ssbull }, d + dF\wedge ).
\tag 3.1.2
$$
Since
$ Y $ is a divisor, we have
$$
\bold{R}\Gamma _{Y}{\Cal O}_{X}[1] = {\Cal H}^{1}\bold{R}\Gamma 
_{Y}{\Cal O}_{X} = j_{*}{\Cal O}_{X\backslash Y}/{\Cal O}_{X},
$$
where
$ j : X \backslash  Y \rightarrow  X $ denotes the inclusion morphism.
Note that a natural morphism
$$
{\Cal H}^{i}p_{+}\bold{R}\Gamma _{Y}{\Cal O}_{X} \rightarrow  
{\Cal H}^{i}p_{+}{\Cal O}_{X}
\tag 3.1.3
$$
is surjective for
$ i \ne  -  n $  and zero otherwise.
Indeed, we have a dense open subvariety
$ U $  of
$ S $  such that the
$ {\Cal H}^{i}p_{+}\bold{R}\Gamma _{Y}{\Cal O}_{X}|_{U} $ are locally free
$ {\Cal O}_{U} $-Modules of finite type, and their fibers at
$ s \in  U $ are isomorphic to
$ {H}_{{Y}_{s}}^{i+n}(X_{s},{\Bbb C}) $.
Then we may replace
$ S $ with
$ U $,
because the
$ {\Cal H}^{i}p_{+}{\Cal O}_{X} $ are simple
$ {\Cal D}_{S} $-Modules (more precisely,
$ {\Cal O}_{S} $ or zero).
But the natural morphism
$ H^{i+n}(X_{s},{\Bbb C}) \rightarrow  
H^{i+n}(X_{s}\backslash Y_{s},{\Bbb C}) $ is an isomorphism for
$ i = -n $,
and zero otherwise.

For the proof of (0.4), we will further simplify the right-hand side of 
(3.1.2).

\proclaim{3.2.~Proposition}
We have natural isomorphisms
$$
\Omega ^{i}({\Bbb A}_{S}^{n+1}/S)^{(m)} = (p\scirc \pi )_{*}{\Omega }
_{{V}^{(m)}/S}^{i}
\tag 3.2.1
$$
compatible with
$ d $,
 $ \wedge  $,
and the action of
$ {\Cal D}_{S} $.
\endproclaim

\demo\nofrills {Proof.\usualspace}
 This is clear if
$ m = 1 $,
because
$ V^{(1)} \rightarrow  {\Bbb A}_{S}^{n+1} $ is the blow up at the origin.
For
$ m > 1 $,
 $ V^{(m)} $ is the quotient of
$ V^{(1)} $ by the natural
$ \mu _{m} $-action on each fiber, where
$ \mu _{m} = \{\lambda  \in  {\Bbb C} : \lambda ^{m} = 1\}. $ So we get 
the assertion taking the invariant part by the action of
$ \mu _{m} $ on
$ {\Bbb A}_{S}^{n+1} $.
\enddemo

\bigskip
\noindent
{\bf 3.3.}
{\it Proof of} (0.4).
By (2.4) applied to
$ M^{\ssbull } = {\Cal O}_{X} $,
 $ we  $ have a canonical morphism of
$ C^{b}({\Cal O}_{X},  \Diff) $
$$
\pi _{*}({\Omega }_{{V}^{(m)}}^{\ssbull }, d + dF\wedge ) \rightarrow  
({\Omega }_{X}^{\ssbull }, d).
\tag 3.3.1
$$
Taking
$ {\Cal H}^{i}p_{+}\DR_{X}^{-1}  $ and using the isomorphism (3.1.2), 
this gives the natural morphism (3.1.3)  by (2.4).
Applying
$ \DR_{S}^{-1} $ (see (1.3.2)) and (1.2.5) to (3.3.1), and then 
transforming right
$ {\Cal D} $-Modules to left
$ {\Cal D} $-Modules, we get a morphism of left
$ p^{-1}{\Cal D}_{S} $-Modules
$$
\pi _{*}({\Omega }_{{V}^{(m)}/S}^{\ssbull }, d + dF\wedge ) \rightarrow  
({\Omega }_{X/S}^{\ssbull }, d).
\tag 3.3.2
$$

Furthermore, the induced morphism
$$
R^{j}(p\scirc \pi )_{*}{\Omega }_{{V}^{(m)}/S}^{i} = R^{j}p_{*}(\pi 
_{*}{\Omega }_{{V}^{(m)}/S}^{i}) \rightarrow  R^{j}p_{*}{\Omega }
_{X/S}^{i}
\tag 3.3.3
$$
is an isomorphism for
$ j \ne  0 $.
Indeed, we have a short exact sequence
$$
0 \rightarrow  \bigoplus _{k\ge 0} {\Omega }_{X/S}^{i}(km) \rightarrow  
\pi _{*}{\Omega }_{{V}^{(m)}/S}^{i} \rightarrow  \bigoplus _{k>0} {\Omega }
_{X/S}^{i-1}(km) \rightarrow  0,
$$
which is the direct image of the short exact sequence
$$
0 \rightarrow  \pi ^{*}{\Omega }_{X/S}^{i} \rightarrow  {\Omega }
_{{V}^{(m)}/S}^{i} \rightarrow  \pi ^{*}{\Omega }_{X/S}^{i-1}\otimes 
_{{\Cal O}_{X}} {\Omega }_{{V}^{(m)}/X}^{1} \rightarrow  0,
$$
because
$ {\Omega }_{{V}^{(m)}/X}^{1} = \pi ^{*}{\Cal O}_{X}(m) $.
So it is enough to show that
$ R^{j}p_{*}{\Omega }_{X/S}^{i}(km) = 0 $ for
$ j, k > 0 $.
But this follows from the Bott vanishing (i.e.,
$ H^{i}({\Bbb P}^{n}, {\Omega }_{{P}^{n}}^{j}(k)) = 0 $ unless  $ i = 0, n $ 
or
$ i = j, k = 0) $,
because
$ R^{n}p_{*}{\Omega }_{X/S}^{i}(km) $ is the dual of
$ R^{0}p_{*}{\Omega }_{X/S}^{n-i}(-km) $ (and
$ R^{0}p_{*}{\Omega }_{X/S}^{n-i} = 0 $ for
$ i \ne n $).

By (3.3.2), we have a morphism between the spectral sequences in the 
category of left
$ {\Cal D}_{S} $-Modules
$$
\align
{E}_{1}^{i,j} = R^{j}(p\pi )_{*}{\Omega }_{{V}^{(m)}/S}^{i} &\Rightarrow  
R^{i+j}(p\pi )_{*}({\Omega }_{{V}^{(m)}/S}^{\ssbull }, d + dF\wedge ),
\tag 3.3.4
\\
{E}_{1}^{i,j} = R^{j}p_{*}{\Omega }_{X/S}^{i} &\Rightarrow  
R^{i+j}p_{*}({\Omega }_{X/S}^{\ssbull }, d),
\tag 3.3.5
\endalign
$$
which are defined by the filtration
$ \sigma  $ on both sides of (3.3.2).
The second spectral sequence (3.3.5) is the Hodge-de Rham spectral 
sequence, and
$ {E}_{1}^{i,j} = 0 $ except for
$ i = j \le  n $ so that it degenerates at
$ E_{1} $.
For (3.3.4) we have
$ {E}_{1}^{i,j} = 0 $ except for
$ j = 0 $ or
$ i = j \le  n $ by the isomorphism (3.3.3) for
$ j > 0 $.
Furthermore we have
$ {E}_{1}^{i,i} = {\Cal O}_{S} \,(1 \le  i \le  n) $ for both.
So the spectral sequence (3.3.4) degenerates at
$ E_{2} $,
because the morphism between
$ {E}_{r}^{i,i} $ should be surjective for any
$ r $ by the surjectivity of (3.1.3) for
$ i > - n $ (since
$ {E}_{r}^{p,2i-p} = 0 $ for
$ p < i) $.
Then
$ {E}_{2}^{i,0} $ of (3.3.4) is identified with the primitive part
$ ({\Cal H}^{i-n}p_{+}\bold{R}\Gamma _{Y}{\Cal O}_{X})^{\prim} $ by 
(3.1.2) and (3.3.1), and we get the assertion.

\bigskip
\noindent
{\bf 3.4.}
{\it Remarks.} (i) If  $ F^{-1}(0) \subset  X $ is (reduced and) smooth 
over
$ S $,
then by Griffiths [5], we have an isomorphism
$$
\Gr_{F}(R^{i}p_{*}{\Omega }_{Y/S}^{\ssbull })^{\prim} = 
{\Cal H}^{i+2}(\Omega ^{\ssbull }({\Bbb A}_{S}^{n+1}/S)^{(m)}, dF\wedge ),
\tag 3.4.1
$$
where
$ \Gr_{F} $  is the direct sum of the graded quotients of the Hodge
filtration.
Furthermore, the image of a vector field
$ \xi  $ on
$ S $ by the Kodaira-Spencer map corresponds to the multiplication by
$ - \xi F $.
Since
$ F^{-1}(0) $ is smooth,
$ \{\partial F/\partial x_{i}\} $ is a regular sequence of
$ {\Bbb C}[x] $,
and it is easy to show
$$
\Gr^{G}{\Cal H}^{i}(\Omega ^{\ssbull }({\Bbb A}_{S}^{n+1}/S)^{(m)}, d + 
dF\wedge ) = {\Cal H}^{i}(\Omega ^{\ssbull }({\Bbb A}_{S}^{n+1}/S)^{(m)}, 
dF\wedge ),
\tag 3.4.2
$$
where
$ G $ is the filtration by degree of
$ x_{1}, \dots , x_{n} $.

(ii) If
$ S = pt $,
let
$ U = F^{-1}(1) \subset  {\Bbb A}^{n+1} $.
Then we have an isomorphism
$$
\tilde{H}^{i}(U,{\Bbb C}) = H^{i+1}(\Omega ^{\ssbull }({\Bbb A}^{n+1}), d + 
dF\wedge )
\tag 3.4.3
$$
by [3].
Since
$ X \backslash  Y $ is the quotient of
$ U $ by the action of
$ \mu _{m} $ which gives also the local monodromy, we get the 
isomorphism
$$
{H}_{Y}^{i}(X,{\Bbb C})^{\prim} = \tilde{H}^{i-1}(X\backslash Y, \bold{C}) 
= H^{i}(\Omega ^{\ssbull }({\Bbb A}_{S}^{n+1}/S)^{(m)}, d + dF\wedge ),
\tag 3.4.4
$$
where the first isomorphism follows from the long exact sequence 
associated with local cohomology.
So (0.4) follows in this case.

(iii) Let
$ \tilde{F} = F + {x}_{n+1}^{m} $,
and define
$ X = {\Bbb P}^{n} = \{x_{n+1} = 0\} \subset  \tilde{X}  = {\Bbb P}^{n+1} $,
and
$$
Y = F^{-1}(0)_{\red} \subset  X,\quad \tilde{Y} = \tilde{F}^{-1}(0)_{\red} 
\subset  \tilde{X},\quad U = F^{-1}(-1) \subset  {\Bbb A}^{n+1},
$$
so that
$ Y = \tilde{Y} \cap  X \subset  \tilde{X} $ and
$ U = \tilde{Y} \backslash  Y $.
Then we have a long exact sequence
$$
\rightarrow  {H}_{Y}^{i}(\tilde{X},{\Bbb C}) \rightarrow  
{H}_{\tilde{Y}}^{i}(\tilde{X},{\Bbb C}) \rightarrow  {H}_{U}^{i}
({\Bbb A}^{n+1},{\Bbb C}) \rightarrow .
$$
which gives
$$
\rightarrow  {H}_{Y}^{i-2}(X,{\Bbb C}) \rightarrow  
{H}_{\tilde{Y}}^{i}(\tilde{X},{\Bbb C}) \rightarrow  H^{i-2}(U,{\Bbb C}) 
\rightarrow .
$$
because
$ {H}_{Y}^{i}(\tilde{X},{\Bbb C}) = {H}_{Y}^{i-2}(X,{\Bbb C}) = H_{2n-i+2}
(Y), {H}_{\tilde{Y}}^{i}(\tilde{X},{\Bbb C}) = H_{2n-i+2}(\tilde{Y}) $,
and
$ {H}_{U}^{i}({\Bbb A}^{n+1},{\Bbb C}) = H^{i-2}(U,{\Bbb C}) $.
By [3], p.~25-26, this induces short exact sequences
$$
0 \rightarrow  {H}_{\tilde{Y}}^{i+2}(\tilde{X},{\Bbb C})^{\prim} 
\rightarrow  \tilde{H}^{i}(U,{\Bbb C}) \rightarrow  {H}_{Y}^{i+1}(X,
{\Bbb C})^{\prim} \rightarrow  0.
\tag 3.4.5
$$
Indeed, we have by (3.4.3-4)
$$
\align
\tilde{H}^{i}(U,{\Bbb C})
&= H^{i+1}(\Omega ^{\ssbull }({\Bbb A}^{n+1}), d + dF\wedge ),
\\
{H}_{Y}^{i+1}(X,{\Bbb C})^{\prim}
&= H^{i+1}(\Omega ^{\ssbull }({\Bbb A}^{n+1})^{(m)}, d + dF\wedge ),
\\
{H}_{\tilde{Y}}^{i+2}(\tilde{X},{\Bbb C})^{\prim}
&= H^{i+2}(\Omega ^{\ssbull }({\Bbb A}^{n+2})^{(m)}, 
d + d\tilde{F}\wedge )
\\
&= \bigoplus_{0<j<m} H^{i+1}(\Omega ^{\ssbull }({\Bbb A}^{n+1})^{(j 
\mod m)}, d + dF\wedge ),
\endalign
$$
where
$ \Omega ^{i}({\Bbb A}^{n+1})^{(j \mod m)} = \sum _{|I|=i} \sum _{|\nu 
|+|I|\equiv j \mod m} {\Bbb C}x^{\nu }dx_{I} $.
Note that
$ (\Omega ^{\ssbull }({\Bbb A}^{n+2}), d + d\tilde{F}\wedge ) $ is the 
Koszul complex for
$ \partial /\partial x_{i} + \partial \tilde{F}/\partial x_{i} 
\,(0 \le  i \le  n + 1) $ so that
$$
(\Omega ^{\ssbull }({\Bbb A}^{n+2}), d + d\tilde{F}\wedge ) = (\Omega 
^{\ssbull }({\Bbb A}^{n+1}), d + dF\wedge ) \otimes _{{\Bbb C}} (\Omega 
^{\ssbull }({\Bbb A}^{1}), d + m{x}_{n+1}^{m-1}dx_{n+1}\wedge ).
$$
(More generally, the Thom-Sebastiani type theorem holds for the 
cohomology of a general fiber of a polynomial map in general [4]).

\bigskip\bigskip
\centerline{{\bf References}}

\bigskip

\item{[1]}
A.~Adolphson and S.~Sperber, Dwork cohomology, de Rham cohomology 
and hypergeometric functions, preprint December 1997.

\item{[2]}
A.~Borel et al., Algebraic $ {\Cal D} $-modules, Perspectives in Math. 2, Academic 
Press, 1987.

\item{[3]}
A.~Dimca, On the Milnor fibration of weighted homogeneous polynomials, 
Compos. Math. 76 (1990), 19--47.

\item{[4]}
A.~Dimca and M.~Saito, On the cohomology of a general fiber of a 
polynomial map, Compos. Math. 85 (1983), 299--309.

\item{[5]}
P.~Griffiths, On the period of certain rational integrals I, II, Ann. Math. 
90 (1969), 460--541.

\item{[6]}
A.~Grothendieck and J. Dieudonn\'e, El\'ements de G\'eom\'etrie 
Alg\'ebrique, Publ. Math. IHES 32 (1967).

\item{[7]}
N.~Katz, On the differential equations satisfied by period matrices, 
Publ. Math. IHES 35 (1968), 223--258.

\item{[8]}
N.~Katz and T.~Oda, On the differentiation of de Rham cohomology 
classes with respect to a parameter, J. Math. Kyoto Univ. 8 (1968), 
199--213.

\item{[9]}
C.~Sabbah, On the comparison theorem for elementary irregular
$ {\Cal D} $-modules, Nagoya Math. J. 141 (1996), 107--127.

\item{[10]}
M.~Saito, Modules de Hodge polarisables, Publ. RIMS, Kyoto Univ., 24 
(1988), 849--995.

\item{[11]}
\SameAuthor, Induced
$ {\Cal D} $-Modules and differential complexes, Bull. Soc. math. France 
117 (1989), 361--387.

\item{[12]}
J.-L.~Verdier,
Cat\'egories d\'eriv\'ees, in SGA 4 1/2, Lect. Notes in Math., Springer, 
Berlin, vol. 569, 1977 , pp. 262--311.

\bigskip

\noindent
Alexandru Dimca and Fay\c cal Maaref

\noindent
Math\'ematiques pures, Universit\'e Bordeaux I

\noindent
33405 Talence Cedex, FRANCE

\noindent
e-mail: dimca\@math.u-bordeaux.fr, maaref\@math.u-bordeaux.fr

\bigskip

\noindent
Claude Sabbah

\noindent
UMR 7640 du CNRS, Centre de Math\'ematiques, Ecole Polytechnique

\noindent
91128 Palaiseau Cedex, FRANCE

\noindent
e-mail: sabbah\@math.polytechnique.fr

\bigskip

\noindent
Morihiko Saito

\noindent
RIMS Kyoto University

\noindent
Kyoto 606--8502 JAPAN

\noindent
e-mail: msaito\@kurims.kyoto-u.ac.jp

\bye